\documentclass{amsart}                   
\usepackage{amsfonts}
\usepackage{amsmath}
\usepackage{amscd}
\usepackage{amssymb}
\usepackage[dvips]{graphicx}
\usepackage[dvips]{epsfig}
\usepackage{amsthm}
\usepackage[all]{xy}

\def\mb#1{\mathbb{#1}}

\newcommand{\Real}{\mathbb{R}}

\newcommand{\Sp}{\mathbb{S}}
\newcommand{\T}{\mathbb{T}}

\def\R#1{\mb{R}^{#1}}

\def\L{\mathcal{L}}
\def\X{\widetilde{X}}

\def\U{\mathcal{U}}

\def\s#1{\mb{S}^{#1}}               

\def\Arrow#1{\overset{#1}\longrightarrow}           

\def\Morf#1#2#3{\xymatrix{  #1 \ar[r]^(.5){#2} & #3}}

\def\im{\text{\rm Im}}              

\def\lift#1{\widetilde{#1}}   
\def\X{\lift{X}}

\def\w{\lift{\omega}}

\def\lm#1#2{#1\times c(#2)}

\newcounter{numero}
\newcommand{\Numero}{\setcounter{numero}{1}(\arabic{numero}) }
\newcommand{\numero}{\addtocounter{numero}{1}(\arabic{numero}) }

\newcounter{letra}
\newcommand{\Letra}{\medskip \setcounter{letra}{1}(\alph{letra}) }
\newcommand{\letra}{\medskip \addtocounter{letra}{1}(\alph{letra}) }

\newcounter{romnumero}

\newcounter{bibnumero}

\theoremstyle{plain}
\newtheorem{Teo}{Theorem}[subsection]
\newtheorem{Coro}[Teo]{Corollary}
\newtheorem{Lema}[Teo]{Lemma}
\newtheorem{Pro}[Teo]{Proposition}

\theoremstyle{definition}
\newtheorem{Def}[Teo]{Definition}

\newtheorem{Ejms}[Teo]{Examples}

\theoremstyle{remark}                               
\newtheorem{Obs}[Teo]{Remark}

\newenvironment{Proof}{ [{\it Proof\/}]\rm\hskip2mm }{\hfill$\square$\vskip2mm}       
\newenvironment{sketch}{ [{\it Sketch of the Proof\/}]\rm\hskip2mm }{\hfill$\square$\vskip2mm}
\def\bdem{\begin{Proof}}
\def\edem{\end{Proof}}
\def\bsketch{\begin{sketch}}
\def\esketch{\end{sketch}}

\begin{document}

\title{On the Functoriality of Stratified Desingularizations}

\author{T. Guardia}
\address{Universidad Central de Venezuela-
Esc. de Matem\'atica.} \email{tomas.guardia@ciens.ucv.ve}

\author{G. Padilla}
\address{Universidad Central de Venezuela-
Esc. de Matem\'atica.} \email{gabriel.padilla@ciens.ucv.ve}

\date{}

\maketitle

\begin{center}
{\small\it To our families.}
\end{center}

\begin{abstract}

This article is devoted to the study of smooth desingularization, which are customary employed in the definition of De Rham Intersection Cohomology with differential forms \cite{illinois}. In this paper we work with the category of Thom-Mather simple spaces \cite{pflaum}, \cite{thom}. We construct a functor which sends each Thom-Mather simple space into a smooth manifold called its primary unfolding. Hence we prove that the primary unfoldings are unique up Thom-Mather isomorphisms.

\end{abstract}

\section*{Introduction}

Stratafied spaces were initially defined by Thom \cite{thom}. For a recent
explanation of the topic see Pflaum \cite{pflaum}. These spaces are
intimately related to the study of intersection homology defined by Goresky and MacPherson \cite{gm1}, which is a loose homology theory since the homotopy axiom does not hold in this context. Nevertheless, it extends the Poincar\'e Duality to
spaces with singular points. We define the intersection cohomology through an unfolding. This is a suitable smooth desingularization
of a pseudomanifolds; it consists in removing the
singular part of a pseudomanifold and gluing two diffeomorphic
copies of the regular part through a hypersurface, for more details
see  Davis \cite{davis},  Ferrarotti \cite{ferrarotti} and Saralegi \cite{illinois}.

In this article we work with simple pseudomanifolds which are just
stratified pseudomanifolds whose depth is $\leq 1$. Simple
pseudomanifolds are almost everywhere smooth spaces, they can be
separated in two disjoint pieces: A dense open subset callet its
regular part, and a closed called its singular part. Their local
model has the form $U\times c(L)$ where $U$ and $L$ are manifolds
and $L$ is compact; the name pseudomanifold arises from this kind of
charts. This article is structured as follows: first we introduce
the notion of simple pseudomanifolds and their unfoldings.  Next we
study the Thom-Mather simple spaces and we prove that each simple
pseudomanifold has an unfolding if and only if it satisfies the
Thom-Mather condition.

In the last section we construct a functor from the category of
Thom-Mather simple spaces to the category of the smooth manifolds.
This is done with the tool of primary unfoldings, which allows us to
work in a\linebreak simpler subcategory of spaces. The unfoldings in
general are neither unique nor funtorial objects; see Brasselet,
Hector and Saralegi \cite{brasselet} and Dalmagro \cite{dalmagro2}
.\,These are the main results of this article: First, every radium
preserving morphism between Thom-Mather simple spaces induces a
smooth map between their respective unfoldings. Second, the primary
unfolding of a Thom-Mather simple space is unique.

All along this article, each time we use the word {\it manifold} we mean a smooth differentiable manifold of class $C^\infty$ without boundary.

\section{The geometry}

\subsection{Simple Pseudomanifolds}
In this section we introduce the family of simple pseudomanifolds;
notice that these are stratified pseudomanifolds whose depth is
$\leq1$. Simple spaces are a previous step in regard of a clear
exposition of pseudomanifolds. For more details see \cite{gm1}
\cite{pflaum}. \begin{Def}\label{def simple spaces}
  A \textbf{simple space} is a 2nd countable metric topological space  $X$ which can be written as the disjoint union of two manifolds $X=R\sqcup\Sigma$, such that
  $R$ is an open dense set (and therefore $\Sigma$ is closed). We refer to $R$ (resp. $\Sigma$) as the \textbf{regular} (resp. \textbf{singular})
    part of $X$.
    A regular (resp. singular) {\bf stratum} $S$ of $X$ is a connected component of $R$ (resp. $\Sigma$). The pair  $(R,\Sigma)$ is a {\bf decomposition} of $X$.\vskip2mm

    A {\bf simple subspace} of $X$ is a subset $Y\subset X$ such that
    $(R\cap Y,\Sigma\cap Y)$ is a decomposition of $Y$ with the induced topology. This condition is not trivial, because according to 1.1.1 $R\cap Y$ and $\Sigma\cap Y$ must be manifolds. \vskip2mm

    If $X'$ is another simple space, then  a {\bf morphism} (resp. \textbf{isomorphism})
  is a continuous map $\Morf{X}{f}{X'}\label{def simple morphism}$ which preserves the decomposition in a smooth (resp. diffeomorphic) way. An \textbf{embedding} is a morphism $\Morf{X}{f}{X'}$
  such that $f(X)$ is a simple subspace of $X'$ and $\xymatrix{X \ar[r]^(.4){f}&  f(X)}$ is an isomorphism.
  \qed
\end{Def}

\begin{Ejms}\label{ejems simple spaces}
    \item[\Numero] Each  manifold $M$ is a simple space
    whose singular part $\Sigma=\emptyset$ is the empty set.
    \item[\numero] The canonical decomposition of any manifold with (nonempty) boundary $M$ is
    $(M-\partial M,\partial  M)$.
    \item[\numero] Any open subspace of a simple space is itself a simple space.
    \item[\numero] If $M$ is a manifold and $X$ is a simple space, then the product
    space $M\times X$ is a simple space, its decomposition is $(M\times R, M\times
    \Sigma)$.
    \item[\numero] Let $L$ be a compact manifold, the
    \textbf{open cone} of $L$ is the quotient space
    \[
        c(L)=\frac{L\times [0,\infty)}{\sim}
    \]
        where $(l,0)\sim (l',0)$ for any $l,l'\in L$.
        The equivalence class of a point $(l,r)$ will be written $[l,r]$.
        The \textbf{vertex} is
        the class of any point $(l,0)$; it will be denoted $v$. For convenience,  we agree
        that $c(\emptyset)=\{v\}$ is a singleton. The space $c(L)$
        is simple, its decomposition is $(L\times\R{+},\{v\})$.
        \item[\numero] A {\bf pseudo-Euclidean model}  (or {\bf pem} for short), is a product
        $U\times c(L)$ with the decomposition given by the above examples, i.e.,
        $(V\times L\times\R{+},U\times\{v\})$; where $U,L$ are manifolds and $L$ is
        compact; this $L$ is said to be the {\bf link} of $U$ on $U\times c(L)$.
        Our convention for $L=\emptyset$ implies that any euclidean
        nbhd is a pem.
        \item[\numero] Since any pem is the quotient of a product manifold; each morphism \linebreak
         $\Morf{\lm{U}{L}}{f}{\lm{V}{N}}$ can be
         written as
        \[
                f(u,[l,r])=(a_1(u,l,r), [a_2(u,l,r), a_3(u,l,r)])
        \]
        where $a_1,a_2,a_3$ are maps defined on $U\times L\times[0,\infty)$, and
        they are "smooth by pieces", i.e., on $U\times L\times \{0\}$ and $U\times L\times \R{+}$. Notice
        that $a_3(u,l,0)=0$ for any $u,l$.
\end{Ejms}

\begin{Def}\label{def spm}
        A {\bf simple pseudomanifold} or {\bf spm} for short, is a simple space $X$
        such that each singular point $x\in \Sigma$ has an open nbhd $x\in V\subset X$ which is the
        image of an embedding
        \[
                \alpha:\Morf{\lm{U}{L}}{}{X}
        \]
        We call $V=\im(\alpha)$ a {\bf pem-nbhd} of $x$, while the pair $(U,\alpha)$ is a \textbf{chart}.
        Since any euclidean nbhd is a pem-nbhd (see example  \S\ref{ejems simple
        spaces}-(6)), the above condition is non-trivial just for singular points.
        \vskip2mm

        Up to some minor
        details, we assume that $\alpha(u,v)=u$ for all $u\in U$, so
        $U=V\cap \Sigma$ is an open nbhd of $x$
        on the corresponding stratum $S\subset\Sigma$ containing $x$. We usually ask the points of $S$
        to have the same link $L$, so it does not depend on the choice of $x$. We call $L$ the {\bf link} of $S$.\qed
\end{Def}

\begin{Ejms}
   \item[\Letra] The examples \S\ref{ejems simple spaces}-(1),(2),(5) and (6) are spm's.
   \item[\letra] If $M$ is a manifold and $X$ is a spm then $M\times X$ is a spm.
   \item[\letra] Any open set of a spm is also a spm.
\end{Ejms}

\section{Desingulatizations}

It is well known how some usual (co)homological properties of smooth
manifolds are lost when we add singularities. This is the case of
the Poincar\'e Duality, for instance. In order to recover these
properties on a larger family of spaces, the original works of
Goresky and MacPherson defined the Intersection Homology with
singular chains \cite{gm1}. Later on \cite{forme diff hector-sar},
\cite{illinois}, Hector and Saralegi provided a smooth approach to
the Intersection Cohomology. For more details see also
\cite{brasselet}, \cite{ferrarotti}. Their viewpoint strongly
depends on two geometric objects associated to any spm $X$, which
are built in order to study the way we reach the singular part and
how we can recover the usual cohomological data; these are the
Thom-Mather tubular nbhds \cite{thom} and the smooth unfoldings
\cite{davis}. In the sequel, we will study both of them.

\subsection{Unfoldings}
Recall  the definition of smooth unfoldings \cite{illinois}. Such an object is obtained from an spm $X$ with decomposition
$(R,\Sigma)$ by gluing a finite number of copies of $R$ and replacing $\Sigma$ with a suitable smooth hypersurface.

\begin{Def}\label{def unfoldings}
     An  \textbf{unfolding} of a spm $X$ is a manifold
    $\widetilde{X}$ together a continuous proper map $\X\Arrow{\mathcal{L}}X$ such that:
    \begin{enumerate}
            \item The restriction $\mathcal{L}^{-1}(R)\Arrow{\mathcal{L}}R$ is a smooth trivial covering.
            \item For each  $z\in \mathcal{L}^{-1}(\Sigma)$ there is a commutative square diagram:
            \[
                \begin{array}{ccc}
            U\times L\times \R{} & \Arrow{\widetilde{\alpha}} &     \X \\
            ^{_c}\downarrow & & \downarrow^{_{\mathcal{L}}}  \\
                U\times c(L) & \Arrow{\alpha} & X
            \end{array}
        \]
        such that
        \begin{enumerate}
            \item $(U,\alpha)$ is a chart (see \S\ref{def spm}).
            \item $c(u,l,t)=(u,[l,|t|])$.
            \item $\widetilde{\alpha}$  is a diffeomorphism  on $\mathcal{L}^{-1}(\text{Im}(\alpha))$.
        \end{enumerate}
    We will refer to the above diagram as an {\bf unfolded chart} at $x=\mathcal{L}(z)$.
    \end{enumerate}
    \vskip2mm
    A spm $X$ is \textbf{unfoldable} if there is (at least) a (smooth) unfolding as above. An
    \textbf{unfoldable morphism} is a commutative square diagram
    \[
        \begin{array}{ccc}
            \X & \Arrow{\widetilde{f}} &    \X' \\
            ^{_{\mathcal{L}}}\downarrow & & \downarrow^{_{\mathcal{L}'}}  \\
                X & \Arrow{f} & X'
            \end{array}
        \]
        such that $f$  is a morphism, $\widetilde{f}$ is smooth
        and the vertical arrows are unfoldings.
\end{Def}

\begin{Ejms}\label{ejems unfoldings}
    \item[\Numero]  For any manifold $M$ the identity map $\Morf{M}{id}{M}$ is an unfolding.
    \item[\numero] The map $c$ given in \S\ref{def unfoldings}-(2).(b) is an unfolding of the pem $U\times c(L)$.
    \item[\numero] For any manifold $M$ with non-empty border $\partial M\neq\emptyset$; the link of $\partial M$ is a point $L=\{l\}$.
    Define $\widetilde{M}=M\times\{\pm 1\}/\sim$ as the quotient set obtained by gluing two copies of $M$ along $\partial M$, i.e.
    $(m,1)\sim(m,-1)\forall m\in\partial M$. Write $[m,j]$ for the class of $(m,j)$. Consider the map
    \[
            \Morf{\widetilde{M}}{\mathcal{L}}{M}\hskip1cm \mathcal{L}([m,j])=m
    \]
    Since the border submanifold $\partial M$ has always a smooth collar then, locally, $\mathcal{L}$ behaves like $c$ at
    \S\ref{def unfoldings}-(2).(b). So $\widetilde{M}$ has a unique smooth structure such that $\mathcal{L}$ is an unfolding.
    \item[\numero] If $\Morf{\X}{\mathcal{L}}{X}$ is an unfolding, then for any open subset $A\subset X$
    the restriction $\xymatrix{\mathcal{L}^{-1}(A)\ar[r]^(.6){\mathcal{L}}& A}$ is an unfolding.
    \item[\numero] If $\Morf{\X}{\mathcal{L}}{X}$ is an unfolding, then for any singular stratum $S\subset\Sigma$
    the restriction $\xymatrix{\mathcal{L}^{-1}(S)\ar[r]^(.6){\mathcal{L}}& S}$ is a smooth bundle with typical fiber $F=L$ and $\mathcal{L}^{-1}(S)$ is a hypersurface of $\X$.
    \item[\numero] The 2-torus  $\T^2$ provides an unfolding for the real projective plane
    $\xymatrix{\T^2 \ar[r]^(.4){\mathcal{L}} & \Real P^2}$.
    We replace the singular point $\Sigma=\{\infty\}$ with two disjoint circles $\Sp^1\times\{\pm1\}$ and glue them
    with two copies of the cylinder $\s1\times\Real$, obtaining so $\T^2$. The map $\mathcal{L}$ is given by
   $\mathcal{L}(e^{i\theta},e^{i\varphi})=[e^{i2\theta},e^{i\varphi}]$.
    \item[\numero] The unfoldings are related to some other geometric objects, such as
    tubular nbhds and normalizations. Although we we do not explicitly mention it, we always can assume that the links are connected
     \cite{normalizer}.
\end{Ejms}

\subsection{Tubular neighborhoods}
The construction  of tubular nbhds is due, among others, to Gleason
and Palais. It was developed in the context of compact
transformation groups. In the smooth context, a tubular can be
constructed trough an invariant riemannian metric and they are
related to the existence of equivariant slices \cite{bredon}. All
these results have been extended to the singular context \cite{forme
diff hector-sar}, \cite{pflaum}. When we move on a spm $X$, as we
approach to the singular part we must preserve the geometric notion
of conical radium. \vskip2mm

\begin{Def}\label{def tubes}
     Given an spm $X$ let us consider a singular stratum $S\subset\Sigma$.
     A \textbf{tubular nbhd} of $S$  is a fiber bundle $\xi=(T,\tau,S,c(L))$ satisfying
    \begin{enumerate}
        \item $T$ is an open nbhd of $S$ in $X$.
        \item The abstract  fiber of $\xi$ is $F=c(L)$, the open cone of the link $L$ of $S$.
        \item $\tau(x)=x$ for any $x\in S$. In other words,
        the inclusion $S \subset T$ is a section of the fiber bundle.
        \item The structure group of $\xi$ is a subgroup of Difeo$(L)$; i. e.,  if $(U,\alpha)$, $(V,\beta)$
        are two bundle charts and $U\cap V\neq\emptyset$ then the change of charts can be written as
        \[
            \xymatrix{\beta^{-1}\alpha:U\cap V\times c(L)\ar[r]^(.5){} &
            U\cap V\times c(L)}\hskip1cm   \beta^{-1}(\alpha(u,[l,r]))=(u,[g_{_{\alpha\beta}}(u)(l),r])
        \]
        where $g_{\alpha\beta}(u)$ is a diffeomorphism of $L$ for all $u\in U\cap V$.
    \end{enumerate}
    The spm $X$  is a \textbf{Thom-Mather} iff every singular
    stratum has a tubular nbhd. If $X$, $X'$ are Thom-Mather spms a \textbf{Thom-Mather morphism} is a morphism $f$ in the sense of \S\ref{def simple morphism} which preserves the tubular nbhd.\label{def morfismo thomather}

\end{Def}
\begin{Obs}\label{obs global radium}
        Given a tubular nbhd $\xymatrix{T\ar[r]^(.5){\tau} & S}$, by \S\ref{def tubes}-(4), there is global sense of
        radium $\xymatrix{T\ar[r]^(.5){\rho} & [0,\infty)}$ in the whole tube $T$ given by $\rho(u,[l,r])=r$, this function
        is the {\bf radium} of $T$. Notice
        that $\rho^{-1}(\{0\})=S$ and $\rho^{-1}(\R{+})=(T-S)$. There is also an action $\xymatrix{\R{+}\times T\ar[r]^(.5){} & T}$
        by {\bf radium stretching}, given by $\lambda\cdot\alpha(u,[l,r])=\alpha(u,[l,\lambda r])$.
\end{Obs}

 \begin{Obs}\label{obs disjoint tubes}
        Given a spm $X$, the singular strata are disjoint connected manifolds. Since $X$ is normal,
        they can be separated trough a disjoint family of open subsets. Therefore, and up to some
        minor details, if $X$ is Thom-Mather then we can find a disjoint family of tubular
        nhbds. This allows us to simplify some things. In the rest of this section we
         fix a Thom-Mather spm $X$ and we will assume, without loss of generality, that $X$ has a
         unique singular stratum $S=\Sigma$, with a tubular nbhd as above.
\end{Obs}

The main goal of this section is to study how the unfoldings and the tubes relate to each
other. We state this relationship in the following
\begin{Teo}\label{teo TM iff unfoldable}
    A spm is Thom-Mather if and only if it is unfoldable.
\end{Teo}

We will prove separatedly each implication, so next we show how to obtain an unfolding from a tubular nbhd. For an equivalent way of constructing unfoldings, the reader can see \cite{davis} \cite{forme diff hector-sar}.

\begin{Pro}\label{pro TM implies unfoldable}
 Every Thom-Mather spm is unfoldable.
\end{Pro}

\bdem
        For the singular stratum $S=\Sigma$ and a the tubular nbhd
        $\Morf{T}{\tau}{S}$, let's fix a bundle atlas
        $\U=\{(U_\alpha,\alpha)\}_{\alpha\in\mathfrak{I}}$.
        We proceed in three steps.\vskip2mm
        \underline{\tt$\bullet$ Unfolding a chart}:
        For any chart $(U,\alpha)\in\U$, the unfolding of $\tau^{-1}(U)$ is the composition
        \begin{equation}\label{eq local chart unfolding}
                \xymatrix{U\times L\times\Real\ar[r]^(.5){c} &U\times c(L) \ar[r]^(.5){\alpha} &
                \tau^{-1}(U)}
        \end{equation}
        where $c$ is the map given at \S\ref{def unfoldings}-(2).(b).
        \vskip2mm

        \underline{\tt $\bullet$ Unfolding the tube:} Define
        \begin{equation}\label{eq manifold of unfolding tube }
                \widetilde{T}=\frac{\underset{\alpha}{\bigsqcup\ }U_\alpha\times L \times \Real}{\sim}
                \hskip2cm
                (u,l,t)\sim(u,g_{\alpha\beta}(u)(l),t)\ \forall\alpha,\beta\hskip2mm \forall u\in
                U_\alpha\cap U_\beta
        \end{equation}
        as the quotient of the disjoint union, with the above equivalence relation.
        Write $[u,l,t]$ for the equivalence class of a triple $(u,l,t)$.
        According to \cite[p. 14]{steenrod}, the above operation defines a fiber bundle
        \begin{equation}\label{eq projection unfolding tube}
                \xymatrix{\widetilde{T} \ar[r]^(.5){\widetilde\tau}& S} \hskip1cm
                \widetilde\tau([u,l,t])=u
        \end{equation}
        with abstract fiber $F=L\times\mathbb{R}$ and the same structure group of $T$.  Since the
        cocycles are
        smooth, notice that $\widetilde{T}$ is a manifold.
        Let's now define
        \begin{equation}\label{eq def unfolding the tube}
                \xymatrix{\widetilde{T}\ar[r]^(.5){\mathcal{L}}&T}\hskip2cm
                \mathcal{L}([u,l,r])=\alpha(u,[l,|t|])\hskip5mm  \forall x\in U_\alpha\ \forall\alpha
        \end{equation}
        In order to show that the above arrow is an unfolding of $T$, the reader only needs to
        check, for any chart $(U_\alpha,\alpha)\in\U$,
        the smoothness of $\widetilde\alpha$ in the  next commutative square diagram:
        \[
                \xymatrix{U_\alpha\times L \times
                \Real\ar[r]^(.6){\widetilde{\alpha}}\ar[d]_{c}
                &\widetilde{T}\ar[d]^{\mathcal{L}}\\
                U_\alpha\times c(L)\ar[r]_(.6){\alpha}& T}
        \]
        where $\widetilde{\alpha}(u,l,r)=[u,l,r]$ is just to pick the respective equivalence
        class.

        \underline{\tt $\bullet$ Unfolding the whole spm $X$:} Remark that
        \[
                \mathcal{L}^{-1}(T-S)=T^{+}\sqcup T^{-}
        \]
        has two connected components, each of them being a smooth bundle over $S$
        with abstract fibers $F^{+}=L\times\R{+}$  and $F^{-}=L\times\R{-}$ respectively.
        These two components are disjoint because the cocycles are radium-independent.
        The unfolding of the whole space $X$ can be made by taking two copies of the regular part
        $R=X-\Sigma$; say $R^{+},R^{-}$, and gluing them together and suitably with
        $\widetilde{T}$ along $\mathcal{L}^{-1}(T-S)$.
\edem

In order to prove the converse statement, let's recall that any manifold $M$ with nonempty boundary $\partial M\neq\emptyset$ has a  \textbf{collar}, i.e. a transverse nbhd given by a smooth embedding $\xymatrix{\Gamma:\partial M\times [0,\infty) \ar[r]&M}$ such that
$\text{Im}(\Gamma)$ is open in $M$ and $\Gamma(m,0)=m$ for all
$m\in\partial M$.

\begin{Pro}\label{pro unfoldable implies TM}
 Every unfoldable pseudomanifold  $X$ is Thom-Mather.
\end{Pro}

\bdem
        Let  $\xymatrix{\widetilde{X}\ar[r]^(.5){\mathcal{L}} & X}$ be an
        unfolding of $X$. Then $\mathcal{L}^{-1}(X-\Sigma)$ is a
        finite trivial smooth covering of $R=(X-\Sigma)$, i.e. a disjoint union of
        finitely many diffeomorphic copies of $R$. Pick one, say $R_0\cong R$, such that
        $\overline{R_0}$ is a manifold with boundary
        $\partial(\overline{R_0})=\mathcal{L}^{-1}(S)$, and take a collar
        $\xymatrix{\mathcal{L}^{-1}(S)\times \Real
        \ar[r]^(.6){\Gamma}&\overline{R_0}}$ of $\mathcal{L}^{-1}(S)$ in $\overline{R_0}$. Define
        $T=\mathcal{L}(\im(\Gamma))$ and
        \[
                \xymatrix{T \ar[r]^(.5){\tau} & S}\hskip1cm \tau(\mathcal{L}(\Gamma(z,r)))=\mathcal{L}(z)
        \]
        Following \cite{forme diff hector-sar}, the above map provides a tubular nbhd of $S$ and
        each unfolded chart as in \S\ref{def unfoldings}-(2) induces a bundle chart
        \[
                \xymatrix{U\times c(L) \ar[r]^(.5){\widehat\alpha} & \tau^{-1}(U)}
                \hskip1cm
                \widehat\alpha(u,[l,r])=\mathcal{L}(\Gamma(\widetilde\alpha(u,l,0),r))
        \]
        We leave the details to the reader.
\edem


\section{Functorial constructions}
In the previous sections we dealed with the existence of unfoldings in terms of the tubular nbhds. Now we will see in more detail their categorical properties. It can be easily deduced from \S\ref{def unfoldings} that the unfoldings are neither unique nor functorial objects.  We will restrict ourselves to a narrower family of spaces in order to develope some ideas concerning the smooth desingularization of pseudomanifolds and their functoriality, when considered as a topological process.

\subsection{Primary unfoldings} Primary unfoldings are the smallest unfoldings one can find for a given pseudoma-nifold.
They where originally presented by Brasselet, Hector and Saralegi \cite{brasselet} and later redefined by Dalmagro
\cite{dalmagro2}, whose points of view constitute the aim of this section.

\begin{Def}\label{def primary unfoldings}
        A \textbf{primary unfolding} is an unfolding
        $\Morf{\X}{\mathcal{L}}{X}$ in our previous sense, such that the preimage of the regular part
           $\mathcal{L}^{-1}(R)=R_0\sqcup R_1$ is a double (smooth, trivial) covering, i.e., the union of
        exactly two diffeomorphic copies of $R$.
\end{Def}

\begin{Obs}\label{obs primary unfoldings are representative}
        The family of primary unfoldings is representative in the category of unfoldable spms.
        Starting from any unfolding $\Morf{\X}{\mathcal{L}}{X}$, we can construct a primary unfolding
        $\Morf{\X'}{\mathcal{L}'}{X}$ by taking a manifold with border $M=\overline{R_0}$ as in the proof
        of \S\ref{pro unfoldable implies TM} and then proceeding as in example
        \S\ref{ejems unfoldings}-(3) in order to get $\Morf{\widetilde{M}}{\mathcal{L}''}{M}$. Then
        $\X'=\widetilde{M}$ and $\mathcal{L}'=\mathcal{L}\mathcal{L}''$ is the composition.
\end{Obs}

\begin{Lema}\label{lema continuous lifting}
        Let $X,X'$ be two unfoldable pseudomanifolds. Then, for any pair of unfoldings
        $\Morf{\X}{\mathcal{L}}{X}$, $\Morf{\X'}{\mathcal{L}'}{X'}$ and any morphism $\Morf{X}{f}{X'}$;
        there is a unique continuous and almost everywhere smooth map
    $\xymatrix{\widetilde{X}\ar[r]^(.5){\widetilde{f}} &\X'}$ such
    that the square diagram
    \[
                \xymatrix{\X\ar[r]^(.5){\widetilde{f}}
                \ar[d]_{\mathcal{L}}& \X'\ar[d]^{\mathcal{L}'}\\X\ar[r]_f&X'}
    \]
        is commutative;\, we call $\widetilde{f}$ the \textbf{lifting}
        of $f$.
\end{Lema}
\bdem
   In order to simplify the exposition, by the above remark \S\ref{obs primary unfoldings are representative}, we assume that
        \[
                \Morf{\X}{\mathcal{L}}{X}\hskip2cm\Morf{\X'}{\mathcal{L}'}{X'}
        \]
        are primary unfoldings.

 \begin{itemize}
                \item[\Letra] \underline{\tt Lifting $f$ on the regular part:}
        If $\Sigma=\emptyset$ then there is no singular part, and $X=R$ is a manifold.
        It follows that $\X=R_0\sqcup R_1$. Take $\widetilde{f}=f\times \text{id}$, then \S\ref{lema continuous lifting}
        trivially holds.  So let's suppose that $\Sigma\neq\emptyset$ and moreover, by \S\ref{obs disjoint tubes}, we will assume
        that $\Sigma=S$ is a single stratum. By these arguments, we have already defined a continuous function $\widetilde{f}$
        satisfying \S\ref{lema continuous lifting} on $\mathcal{L}^{-1}(R)$. Therefore, we only must find a continuous extension
        of $\widetilde{f}$ to the entire $\X$, i.e., to $\mathcal{L}^{-1}(S)$.

                \item[\letra] \underline{\tt Extension of the lifting:}
                Pick some $\widetilde{z}\in\L^{-1}(S)$. We must show a way to choose $\widetilde{f}(\widetilde{z})$.
                For this sake, let  $\{\widetilde{z}_n\}_n\subset\L^{-1}(R)$ be a sequence converging to $\widetilde{z}$.
                Since $\mathcal{L},f$ are continuous maps and $\mathcal{L}'$ is a continous and proper map; by an argument of compactness,
                and up to some little adjusts, we may assume that the sequence
                $\{\widetilde{f}(\widetilde{z}_n)\}_n$
                converges in $\X'$. We define
                \[
                    \widetilde{f}(\widetilde{z})=\lim\limits_{n\rightarrow\infty} \widetilde{f}(\widetilde{z}_n)
                \]
                If our limit-definition makes sense then it is also continuous; so next we will show the non ambiguity of $\widetilde{f}$.
                Since the former is a local definition, we first study the
                \item[\letra]\underline{\tt Local writing of the lifting:}
                By \S\ref{def unfoldings}-(2) we can
                restrict to unfoldables charts; so we will assume that X$=U\times c(L)$ and $Y=V\times c(N)$
                are trivial pem nbhds and their respective unfoldings are the canonical ones - see \S\ref{ejems unfoldings}-(2).
                Then $f$ can be written as in \S\ref{ejems simple spaces}-(7). The point
                \[
                        \widetilde{z}=(u,l,0)\in U\times L\times\{0\}
                \]
                is the limit of a sequence
                \[
                        \{\widetilde{z}_n=(u_n,l_n,t_n)\}_n\subset U\times L\times(\mathbb{R}-\{0\})
                \]
                So the sequences $\{u_n\}_n$, $\{l_n\}_n$ and $\{t_n\}_n$
                respectively converge to $u$, $l$ and 0. Since $\mathcal{L}, f$ are continuous maps, the sequence
                \[
                        w_n=f(\mathcal{L}(\widetilde{z}_n))=(a_1(u_n,l_n,|t_n|), [a_2(u_n,l_n,|t_n|),a_3(u_n,l_n,|t_n|)])
                \]
                converges to $w=f(\mathcal{L}(\widetilde{z}))=(a_1(u,l,0), v)$. By the continuity of the functions $a_j$
                for $j=1,2,3$ and up to some little adjust on $a_2$ concerning the compactness arguments;
                we get that
                \[
                        \widetilde{w}_n=(a_1(u_n,l_n,|t_n|), a_2(u_n,l_n,|t_n|),\pm a_3(u_n,l_n,|t_n|))
                \]
                converges to $\w=(a_1(u,l,0), a_2(u,l,0),0)$.
                \item[\letra] \underline{\tt The lifting is well defined:}
                From the continuity of the functions $a_i$,
                the element $\widetilde{w}$ does not depend on the choice of a particular sequence $\{\widetilde{z}_n\}_n$.
        \end{itemize}

        Notice that, the lifting $\widetilde{f}$ is always smooth on
    $\mathcal{L}^{-1}(X-\Sigma)$, the preimage of the regular part.

\edem

\begin{Obs}\label{obs bubble labeling}
  Given an unfolding $\Morf{\X}{\mathcal{L}}{X}$; a \textbf{bubble} of $\widetilde{X}$ is a
  connected component of $\mathcal{L}^{-1}(X-\Sigma)$.The above proof is still
    valid if we take any other permutation of the bubbles.
   Along the rest of this paper, we assume that we are working with the identity permutation, unless
   we state the opposite.
\end{Obs}
\begin{Def}
A morphism $f$ between pems nhbds is \textbf{liftable} if its
lifting $\widetilde{f}$ is globally smooth on every $\widetilde{X}$.
This is equivalent to ask $\widetilde{f}$ to be smooth\, on a
nbhd of $\mathcal{L}^{-1}(\Sigma)$.
\end{Def}
\begin{Pro}\label{Prop local smooth lifting requirements}
            A morphism between pem nhbds
        \[
                \Morf{U\times c(L)}{f}{U'\times c(L')}\hskip1cm f(u,[l,r])=(a_1(u,l,r), [a_2(u,l,r), a_3(u,l,r)])
        \]
        is liftable into
        \[
                \Morf{U\times L\times\mathbb{R}}{\widetilde{f}}{U'\times L'\times\mathbb{R}}\hskip1cm
                \widetilde{f}(u,l,t)=(\widetilde{a}_1(u,l,t), \widetilde{a}_2(u,l,t), \widetilde{a}_3(u,l,t))
        \]
        if and only if
        \begin{enumerate}
                \item[\Letra] $\widetilde{a}_1$, $\widetilde{a}_2$ are smooth even extensions of, respectively, $a_1$, $a_2$.
                \item[\letra]  $\widetilde{a}_3$ is either an odd (and therefore smoooth) extension of $a_3$ or it is a smooth even extension
                and $\widetilde{a}_3(u,l,0)=a_3(u,l,0)=0$ for all $u,l$.
        \end{enumerate}
\end{Pro}

\bdem
        If $\widetilde{f}$ is a lifting of $f$, then
        $fc=c'\widetilde{f}$ where $c$ and $c'$ are canonical unfoldings as  in \S\ref{def unfoldings}-(2).
        Checking both sides of this equality we get
        \[
                f(c(u,l,t))=f(u,[l,|t|])=(a_1(u,l,|t|), [a_2(u,l,|t|), a_3(u,l,|t|)])
        \]
        and
        \[
                c'(\widetilde{f}(u,l,t))=c'(\widetilde{a}_1(u,l,t), \widetilde{a}_2(u,l,t), \widetilde{a}_3(u,l,t))=
                (\widetilde{a}_1(u,l,t), [\widetilde{a}_2(u,l,t), |\widetilde{a}_3(u,l,t)|])
        \]
        we conclude that
        \[
                (a_1(u,l,|t|), [a_2(u,l,|t|), a_3(u,l,|t|)])=(\widetilde{a}_1(u,l,t), [\widetilde{a}_2(u,l,t), |\widetilde{a}_3(u,l,t)|])
        \]
        There are two cases; $t=0$ and $t\neq 0$, from which we get \S\ref{Prop local smooth lifting requirements}.
\edem

\smallskip

\begin{Lema}\label{pro liftable cocycles}
Let $X$ be a Thom-Mather spm, then the cocycles of a tubular
neighborhoood are liftable.
\end{Lema}

\bdem
      It is enough to take $\varphi=\beta^{-1}\alpha$ as in \S\ref{def
      tubes}-(4) and check that $a_1(u,l,r)=u$,
      $a_2(u,l,r)=g(u)(l)$,
       $a_3(u,l,r)=r$ satisfy the hypothesis of \S\ref{Prop local smooth lifting
      requirements}.

\edem

\begin{Lema}\label{def conmutativity of the lifting with cocycles}
Let $\xymatrix{f,f':U\times c(L)\ar[r]&U'\times c(L')}$ be
morphisms,\linebreak $\xymatrix{\varphi:U\times c(L)\ar[r]& U\times
c(L)}$, $\xymatrix{\varphi':U'\times c(L')\ar[r]&U'\times c(L')}$ be
isomophisms as in \S\ref{def tubes}-(4). Then $f'\varphi=\varphi'f$
if and only if $a_1,a_3$ are invarant respect to the group action on
the coordinate $l$ and $a_2$ conmutes with the cocycles i.e,
$a_1,a_2$ and $a_3$ satisfy:
\begin{equation}\label{eq cocycle conditions of the lifting}
\begin{split}
  a_1(u,l,r)=a_1'(u,g(u)(l),r)\\
  g'(a_1(u,l,r))a_2(u,l,r)=a_2'(u,g(u)(l),r)\\
  a_3(u,l,r)=a_3'(u,g(u)(l),r)
\end{split}
\end{equation}

\end{Lema}
\bdem
 Let $f$ as above, and $\varphi,\varphi'$ be cocycles as in \S\ref{def tubes}-(4), then

\[\begin{array}{rl}
f'(\varphi(u,[l,r]))&=f'(u,[g(u)(l),r])\\&=(a_1'(u,g(u)(l),r),[a_2'(u,g(u)(l),r),a_3'(u,g(u)(l),r)])
\end{array}\] On the another hand
\[\begin{array}{rl}
\varphi'(f(u,[l,r]))&=\varphi'(a_1(u,l,r),[a_2(u,l,r),a_3(u,l,r)])\\&=
(a_1(u,l,r),[g'(a_1(u,l,r)(a_2(u,l,r)),a_3(u,l,r))])
\end{array}
\]

If $f'\varphi=\varphi'f$ then, after checking the cases $r=0$ and $r\neq 0$, we get \S\ref{def conmutativity of the lifting with
cocycles}.

\edem

\begin{Lema}\label{def global thomather implies local thomather}
Let $\Morf{X}{\Psi}{X'}$ be a Thom-Mather morphism (see \S$\ref{def morfismo thomather}$) then the
restriction of $\Psi$ on each respective local trivialization
satifies \S\ref{def conmutativity of the lifting with cocycles}.

\end{Lema}\label{teo thomather maps are liftable}
\bdem Without any loss of generality, we can assume that $T=X$ and
$T'=X'$ (see $\S\ref{obs disjoint tubes}$). Let $\alpha,\beta$ and
$\alpha',\beta'$ be two bundle charts of $X$ and $X'$ respectively
defined on $U\times c(L)$ and $U'\times c(L')$.

If $\varphi=\beta^{-1}\alpha$ and $\varphi'=(\beta')^{-1}\alpha'$
are their respective cocycles (c.f., \S\ref{def tubes}-(4)). It
follows that $f=(\alpha')^{-1}\Psi\alpha$ and
$f'=(\beta')^{-1}\Psi\beta$ satisfy $\ref{def conmutativity of the
lifting with cocycles}$. \edem

These results imply that there is a functor from the category of Thom-Mather spaces
to the category of the smooth manifolds. Formally we state the next theorem.

\begin{Teo}\label{teo thommather morphism lift into smooth map}
Every Thom-Mather morphism has a smooth Thom-Mather lifting.
\end{Teo}
\bdem
  Let $X, X'$ be two Thom-Mather spms, if $f$ has a smooth lifting $\widetilde{f}$ as in
  \S\ref{Prop local smooth lifting requirements} then, according to $\S\ref{lema continuous
  lifting}$ there is a unique almost everywhere smooth map $\Morf{\widetilde{X}}{\widetilde{f}}
  {\widetilde{X'}}$. by $\S\ref{teo thomather maps are liftable}$
  $\widetilde{f}$ is locally smooth on a nbhd of
  $\mathcal{L}^{-1}(\Sigma)$ then $\widetilde{f}$ is smooth.
\edem

\subsection{Uniqueness of the primary unfoldings}

Theorem \S\ref{teo thommather morphism lift into smooth map} shows that the primary unfoldings have a quite nice,
functorial behaviour.

\begin{Teo}\label{lema isomorphism lifts in diffeomorphism}
  Let $X,X'$ be two Thom-Mather spaces and let $\Morf{X}{f}{X'}$ be
  a Thom-Mather isomorphism then the smooth lifting of $f$,
  $\widetilde{f}$ is a diffeomorphism between the manifolds
  $\widetilde{X}$ and $\widetilde{X'}$.
\end{Teo}

\bdem
 We know that by $\S\ref{lema continuous lifting}$ the lifting of
 $f$, $\widetilde{f}$ is unique up to permutation of the bubbles. Let us fix a bubble permutation $\sigma$( as in \S\ref{obs bubble labeling}) on
 $\mathcal{L}^{-1}(X-\Sigma)$. As $f$ is a Thom-Mather morphism its
 lifting $\widetilde{f}$ is a smooth Thom-Mather map between
 $\widetilde{X}$  and $\widetilde{X'}$ (see, \S\ref{teo thommather morphism lift into smooth
 map}). On the another hand, the inverse morphism of $f$, $g$ is
 also a Thom-Mather map which lifts into a smooth manifold map
 $\widetilde{g}$ between $\widetilde{X'}$ and $\widetilde{X}$. In order to satisfy the equations $\widetilde{g}\widetilde{f}=id_{\widetilde{X}}$
 and $\widetilde{f}\widetilde{g}=id_{\widetilde{X'}}$ we take the lifting $\widetilde{g}$ of $g$ induced with the inverse permutation $\sigma^{-1}$ of the bubbles. So, therefore $\widetilde{f}$ is a
 diffeomorphism.
\edem

\begin{Coro}
The primary unfolding of a Thom-Mather space is unique.
\end{Coro}

\bdem Take $X=X'$ as in $\S\ref{lema isomorphism lifts in
diffeomorphism}$ and apply the same argumentation of this lemma to
the identity map on $X$.According to \S\ref{obs bubble labeling}, this is true up to bubble labeling.

 \edem

\section{Acknowledgments}

T. Guardia thanks the CDCH-UCV and the FONACIT, so as G. Padilla
thanks the Centro de Geometr\'{\i}a UCV and the Labo CGGA;
institutions which have partially supported this research. Both
authors would also like to thank the hospitality of the staff in the
IVIC's Biblioteca Marcel Roche.

\end{document}